\documentclass{article} 
\usepackage{graphicx}
\usepackage{amsmath}
\usepackage{amsfonts}
\usepackage{amssymb}
\usepackage{latexsym}
\usepackage{tikz}
\usepackage{fancyvrb}
\usepackage{datetime}

\usepackage{color} 
\usepackage[colorlinks]{hyperref}
\usepackage{cite}

\newtheorem{lemma}{Lemma}[section]
\newtheorem{proposition}{Proposition}[section]

\newtheorem{remark}{Remark}
\newtheorem{theorem}{Theorem}[section]
\newtheorem{conjecture}{Conjecture}

\usepackage{graphicx}
\usepackage{amsmath}
\usepackage{amsfonts}
\usepackage{amssymb}
\usepackage{latexsym}

\begin{document}
 \begin{center}
{\Large{The spectral determinations of the multicone graphs $ K_w\bigtriangledown mC_n$}

\medskip 

{Ali Zeydi Abdian\footnote{Lorestan University, College of Science, Lorestan, Khoramabad, Iran; e-mail: aabdian67@gmail.com; azeydiabdi@gmail.com}}
} \end{center}
 \begin{abstract}
 The main goal of the paper is to characterize new classes of multicone graphs which are determined by both adjacency and Laplacian spectra. A multicone graph is defined to be the join of a clique and a regular graph.  A wheel graph obtained from the join of a complete graph on a vertex with a cycle. A question about when wheel graphs are determined by their adjacency spectra is still unsolved. So, any indication about the determinations of  these graphs with respect to their adjacency spectra can be an interesting and important problem. In [{\it Y. Zhang, X. Liu, and X. Yong:}  Which wheel graphs are determined by their Laplacian spectra?. Comput. Math. Appl., {\bf 58} (2009) 1887--1890] and  [{\it M.-H. Liu:} Some graphs determined by their (signless) Laplacian spectra. Czech. Math. J., {\bf 62}, (2012) 1117--1134]  it have been shown that except for, the wheel graph of  order seven, all wheel graphs are determined by their Laplacian spectra and wheel graphs are determined by their signless Laplacian spectra, respectively. In this study, we present new classes of connected multicone graphs which are a natural generalization of wheel graphs and we show that these graphs are determined by their adjacency spectra as well as their Laplacian spectra.   Also, we show that  complement of some of these graphs are determined by  their adjacency spectra. In addition, we give a necessary and sufficient condition for perfecting graphs cospectral with presented graphs in the paper. Finally, we pose two problems for further work.
\end{abstract}

\section{ Introduction} 

In this paper, all graphs, except in Section 5, are  connected undirected simple graphs (for answering to this question why we consider connected graphs, one can see the paragrah before Theorem \ref{the 2-2} and also Theorem \ref{the 2-2} of the paper). Let $G = (V (G),E(G))$ be a graph with vertex set $ V = V (G) = \left\{ {v_1, . . . , v_n} \right\}$ and edge set $E(G)$. All graphs considered here are simple and undirected. All notions on graphs that are not defined here can be found in \cite{Ba, B, CRS, K, W}. A graph consisting of $ k $ disjoint copies of an arbitrary graph $ G $ will be denoted by the $ kG $. The complement of a graph $ G $ is denoted by $ \overline{G} $. The join of two graphs $ G $ and $ H $ is the graph obtained  from disjoint union of $ G $ and $ H $ and connecting any vertex of $ G $ to any vertex of $ H $. The join of two graphs $ G $ and $ H $ is denoted by $ G\bigtriangledown H $. We say that a graph $ G$ is an $ r $-regular graph, if the degree of its regularity is $ r $. Given a graph $ G $, the cone over $ G $ is the graph formed by adjoining a vertex adjacent to every vertex of $ G $.
Let the matrix $A(G)$ be the $(0,1)$-adjacency matrix of $ G$ and $d_k$ be the degree of the vertex $v_k$. The matrix $L(G) = D(G)-A(G)$ is called the Laplacian matrix of $G$, where $D(G)$ is the $ n\times n $ diagonal matrix with $ V = V (G) = \left\{ {d_1, . . . , d_n} \right\}$ as diagonal entries (and all other entries $ 0 $). Since both matrices $A(G)$ and $L(G)$ are real and symmetric, their eigenvalues are all real numbers. Assume that $\lambda_1\geq \lambda_2 \geq . . . \geq \lambda_n $ and $\mu_1\geq \mu_2\geq . . . \geq \mu_n (= 0) $ are respectively
the adjacency eigenvalues and the Laplacian eigenvalues of graph $ G $. The adjacency spectrum of the graph $ G $ consists of the adjacency eigenvalues (together with their multiplicities), and the Laplacian spectrum of the graph $ G $ consists of the Laplacian eigenvalues (together with their multiplicities) and we denote them by $ \operatorname{Spec}_{A}(G) $ and $ \operatorname{Spec}_{L}(G) $, respectively. Two graphs $ G $ and $ H $ are said to be cospectral if they have equal spectrum (i.e.,
equal characteristic polynomial).  If $ G $ and $ H $ are isomorphic, they are necessarily cospectral. Clearly, if two graphs are cospectral, they must possess equal number of vertices. We say that a graph $ G $ is determined by its adjacency (Laplacian) spectra (DS, for short), if for any graph $ H $ with $ \operatorname{Spec}_{A}(G)=\operatorname{Spec}_{A}(H)$ ($ \operatorname{Spec}_{L}(G)=\operatorname{Spec}_{L}(H)$), $ G$ is isomorphic to $ H $. \\
So far numerous examples of cospectral but non-isomorphic
graphs are constructed by interesting techniques such as Seidel switching, Godsil-McKay switching, Sunada or Schwenk method. For more information, one may see \cite{VH, VH1} and the references cited in them. Only a few graphs with very special structures have been reported to be determined by their spectra (DS, for short) (see \cite{A, AA, AAA, AAAA, BJ, BH, CHVW, D, HLZ, P1, LS, M, VH, VH1, WBHB, W2} and the references cited in them). Recently Wei Wang and Cheng-Xian Xu have developed a new method in \cite{VH} to show that many graphs are determined by their spectrum and the spectrum of their complement.
Van Dam and Haemers \cite{VH} conjectured that almost all graphs are determined by their spectra. Nevertheless, the set of graphs that are known to be determined by their spectra is too small. So, discovering classes of graphs that are determined by their spectra can be an interesting problem. The characterization of DS graphs goes back about half of a century and it originated in Chemistry \cite{P1, P}. About the background of the question ''Which graphs are determined by their spectrum?'', we refer to \cite{VH}. A spectral characterization of multicone graphs were studied in \cite{W1, WBHB}. In \cite{W1}, Wang, Zhao and Huang investigated on the spectral characterization of multicone graphs and also they claimed that friendship graphs $ F_n $ (that are special classes of multicone graphs) are DS with respect to their adjacency spectra. In addition, Wang, Belardo, Huang and Borovi\'canin \cite{WBHB} proposed such conjecture on the adjacency spectrum of $ F_n $. This conjecture caused some activities on the spectral characterization of $ F_n $.  Finally, Cioab\"a and et al., \cite{CHVW} proved that if $ n\neq 16 $, then friendship graphs $ F_n $ are DS with respect to their adjacency spectra. Abdian and Mirafzal \cite{A} characterized new classes of multicone graphs which were DS with respect to their spectra. Abdian \cite{AA} characterized two classes of multicone graphs and proved that the join of an arbitrary complete graph and the generalized quadrangle graph $ GQ(2,1) $ or $ GQ(2,2) $ is DS with respect to its adjacency spectra as well as its Laplacian spectra. This author also proposed four conjectures about adjacency spectrum of complement and signless Laplacian spectrum of these multicone graphs. In \cite{AAA}, the author showed that multicone graphs $ K_w\bigtriangledown P_{17} $ and $ K_w\bigtriangledown S $ are DS with respect to their adjacency spectra as well as their Laplacian spectra, where $ P_{17} $ and $ S $ denote the Paley graph of order 17 and the Schl\"afli graph, respectively. Also, this author conjectured that these multicone graphs are DS with respect to their signless Laplacian spectra. In \cite{AAAA}, the author proved that multicone graphs $ K_w\bigtriangledown L(P) $ are DS with respect to both their adjacency and Laplacian spectra, where $ L(P) $ denotes the line graph of the Petersen graph. He also proposed three conjectures about the signless Laplacian spectrum and the complement spectrum of these multicone graphs. For getting further information about characterizing some multicone graphs which are DS see \cite{AAAAA, AAAAAA}.\\ 
We believe that the proofs in \cite{W1} contain some gaps. In \cite{W1}, the authors conjectured that if a graph is cospectral to a friendship graph, then its minimum degree is $ 2 $ (see Conjecture 1). In other words, they could not determine the minimum degree of graphs cospectral to a (bidegreed) multicone graph (see Conjecture 1). Hence, by their techniques (\cite{W1}) cannot characterize new classes of multicone graphs that we want to characterize them. Conjectures (Conjectures 1 and 2) which had been proposed by Wang, Zhao and Huang \cite{W1} are not true and there is a counterexample for them (see the first paragraph after Corollary 2 of \cite{CHVW}). In Theorem 3 (ii) of \cite{W1} first the minimum degree of a graph cospectral to a graphs belonging to $ \beta(n-1, \delta) $ (classes of bidegreed graphs with degree sequence $ \delta $ and $ n-1 $, where $ n $ denotes the number of vertices) must be determined, since in general the minimum degree of a graph cannot be determined by its spectrum. Therefore, we think that theorem without knowing the minimum degree of a graph cospectral with one of graphs $ \beta(n-1, \delta) $ will not be effective and useful.\\
In this paper, we present some techniques which enable us to characterize graphs 
 that are DS with respect to their adjacency and Laplacian spectra.\\

This paper is organized as follows. In Section 2, we review some basic information and preliminaries. In Section 3, we present new classes graphs that are either bidegreed or regular  with respect to their adjacency spectrum. In Section 4, we prove that any  graph cospectral with one of these graphs are determined by their adjacency spectrum.  In Section 5, we show that complement of some classes of these graphs are determined by their adjacency spectrum. In Section 6, we prove that these graphs are DS with respect to their Laplacian spectrum. In Section 7, we show that any graph cospectral with special classes of these graphs must be perfect. In Section 8, we review what were said  in the previous sections and finally we propose two conjectures for further research.

\section{ Preliminaries}

In this section we present some results which will play an important role throughout this paper.

\begin{lemma}\cite{A, AA, AAA, AAAA, M, VH}\label{lem 2-1}
Let $ G $ be a graph. For the adjacency matrix and the  Laplacian matrix of $ G $, the following can be obtained from the spectrum: 

\begin{flushleft}
$(i)$ The number of vertices,

$(ii)$ The number of edges.

For the adjacency matrix, the following follows from the spectrum:

$(iii)$ The number of closed walks of any length,

$(iv)$ Whether $G$ is regular, and the common degree,

$(v)$ Being bipartite or not.

For the Laplacian matrix, the following follows from the spectrum:

$(vi)$ The number of spanning trees,

$(vii)$ The number of components,

$(viii)$ The sum of squares of degrees of vertices.
\end{flushleft}
\end{lemma} 

\begin{theorem}\cite{A, AA, AAA, AAAA, CRS, M, W1}\label{the 2-100} If $G_1$ is $r_1$-regular with $n_1$ vertices, and $G_2$ is $r_2$-regular with $n_2$ vertices, then the characteristic polynomial of the join $ G_1\bigtriangledown G_2 $ is given by:

\begin{center}
$P_{G_{1}\bigtriangledown G _{2}(y)}=\frac{P_{G_{1}}(y)P_{G_{2}}(y)}{(y-r_1)(y-r_2)}((y-r_1)(y-r_2)-n_1n_2)$.

\end{center}

\end{theorem}

The spectral radius of a graph $ \Lambda $ is the largest eigenvalue of adjacency matrix of graph $ \Lambda $ and it is denoted by $ \varrho(\Lambda) $. A graph is called bidegreed, if the set of degrees of its vertices consists of two elements.\vspace{2mm}

For further information about the following inequality we refer the reader to  \cite{W1}  (see the first paragraph after Corollary 2.2 and also Theorem \label{the 2-1} of \cite{W1}). It is stated in \cite{W1} that if $G$ is disconnected, then the equality in the following can also occur. However, in this paper we only consider {\bf connected case} and we state the equality in this case.
\begin{theorem}\cite{A, AA, AAA, AAAA,  M, W1}\label{the 2-2} 
Let $ G $ be a simple graph with $ n$ vertices and $ m $ edges. Let $ \delta=\delta(G)$ be the minimum degree of vertices of $ G $ and $\varrho(G)$ be the spectral radius of the adjacency matrix of $ G $. Then 
\begin{center}
$ \varrho(G)\leq \frac{\delta-1}{2}+\sqrt{2m-n\delta+\frac{(\delta+1)^{2}}{4}} $.

\end{center}
Equality holds if and only if $ G $ is either a regular graph or a bidegreed graph in which each vertex is of degree either $ \delta $ or $ n-1 $.
\end{theorem}

\begin{theorem} \cite{A, AA, AAA, AAAA,  M}\label{the 2-3} Let $ G $ and $ H $ be two graphs with Laplacian spectrum $ \lambda_1\ge\lambda_2\ge...\ge\lambda_n $ and $ \mu_1\ge\mu_2\ge...\ge\mu_m $, respectively. Then Laplacian spectra of $ \overline{G}\ $ and $G \bigtriangledown H $ are $ n -\mathop \lambda \nolimits_1 ,n - \mathop \lambda \nolimits_2 ,...,n - \mathop \lambda \nolimits_{n - 1} ,0 $ and $ n + m,m + \mathop \lambda \nolimits_1 ,...,m + \mathop \lambda \nolimits_{n - 1} ,n + \mathop \mu \nolimits_1 ,. . ., n + \mathop \mu \nolimits_{m - 1} ,0$, respectively.
\end{theorem}

\begin{theorem}\cite{A, AA, AAA, AAAA,  M}\label{the 2-4} Let $ G $ be a graph on $n$ vertices. Then $n$ is one of the Laplacian eigenvalue of $ G $ if and only if $ G $ is the join of two graphs. 
\end{theorem}

\begin{theorem}\cite{K}\label{the 2-5} For a graph $G$, the following statements are equivalent: 

\begin{flushleft}
$(i)$ $G$ is $d$-regular.

$(ii)$ $ \varrho(G)=d_{G} $, the average vertex degree.

$(iii)$ $G$ has $ v=(1,1,...,1)^T $ as an eigenvector for $ \varrho(G)$.

\end{flushleft}
\end{theorem}

\begin{proposition}\cite{A, AA, AAA, AAAA, CRS,  M, R}\label{prop 2-1} Let $ G-j $ be the graph obtained from $ G $ by deleting the vertex $ j $ and all edges containing $ j $. Then $ P_{G-j}(y)=P_{G}(y)\sum\limits_{i = 1}^m {\frac{{\alpha^2 _{ij}}}{{y - {\mu _i}}}} $, where $ m $ and $ \alpha _{ij} $ are the number of distinct eigenvalues and the main angles (see \cite{R}) of graph $ G $, respectively. 
\end{proposition}
\begin{proposition}\label{prop 2-2}\cite{VH1} Let $ G $ be a disconnected graph that is determined by the Laplacian spectrum. Then the cone over $ G $, the graph $ H $; that is, obtained from $ G $ by adding one vertex that is adjacent to all vertices of $ G $, is also determined by its Laplacian spectrum. 
\end{proposition}

\begin{theorem}\cite{A}\label{the 2-6}
Let $ \Gamma $ be a non-regular graph with three distinct eigenvalues
$ \theta_0>\theta_1>\theta_2$. Then the following hold:

\begin{flushleft}
$(i)$ $ \Gamma $ has diameter two.

$(ii)$ If $ \theta_0 $ is not an integer, then $ \Gamma $ is complete bipartite.

$(iii)$ $ \theta_1\ge 0$ with equality if and only if $ \Gamma $ is complete bipartite.

$(iv)$ $\theta_2 \leq -\sqrt{2}$ with equality if and only if $ \Gamma $ is the path of length $ 2$.
\end{flushleft}
\end{theorem}
\begin{theorem}\cite{A}\label{the 2-7}
 A graph has exactly one positive eigenvalue if and only if its non-isolated
vertices form a complete multipartite graph.
 \end{theorem}

\begin{flushleft}
\begin{remark}
In the following, we always suppose that $w$ and  $n\geq 3$  are natural numbers. Also, $ C_n $ and $ K_w $ denote a cycle of order $ n $ and a complete graph on $ w $ vertices, respectively.  In addition, when ${\rm{Spec}}_{A}(G)={\rm{Spec}}_{A}(K_w\bigtriangledown mC_n)  $,  we always suppose that $ G $ is connected. Because, there are some classes of disconnected graphs  that are not  detemined by their spectra. For example,   $ {\rm{Spec}}_{A}((2C_4\bigtriangledown (3C_4\cup K_3))\cup 5C_4)={\rm{Spec}}_{A}(K_3\bigtriangledown 10C_4) $ but $ (2C_4\bigtriangledown (3C_4\cup K_3))\cup 5C_4\ncong K_3\bigtriangledown 10C_4 $, ${\rm{Spec}}_{A}((C_5\bigtriangledown (6C_5\cup K_3))\cup 4C_5)={\rm{Spec}}_{A}(K_3\bigtriangledown 11C_5) $ but $ C_5\bigtriangledown (6C_5\cup K_3))\cup 4C_5 \ncong K_3\bigtriangledown 11C_5$ and ${\rm{Spec}}_{A}((C_6\bigtriangledown (2C_6\cup K_3))\cup 2C_6)={\rm{Spec}}_{A}(K_3\bigtriangledown 5C_6) $ but $ C_6\bigtriangledown (2C_6\cup K_3))\cup 2C_6 \ncong K_3\bigtriangledown 5C_6$. If graphs cospectral with one of multicone graphs $ C_3\bigtriangledown mC_n $ are connected, then these  graphs are determined by their adjacency spectrum (see Theorem 4.1).
\end{remark} 
\end{flushleft}

\hspace{10mm} 
\begin{figure}[h] 
\centerline{\includegraphics[height=5cm]{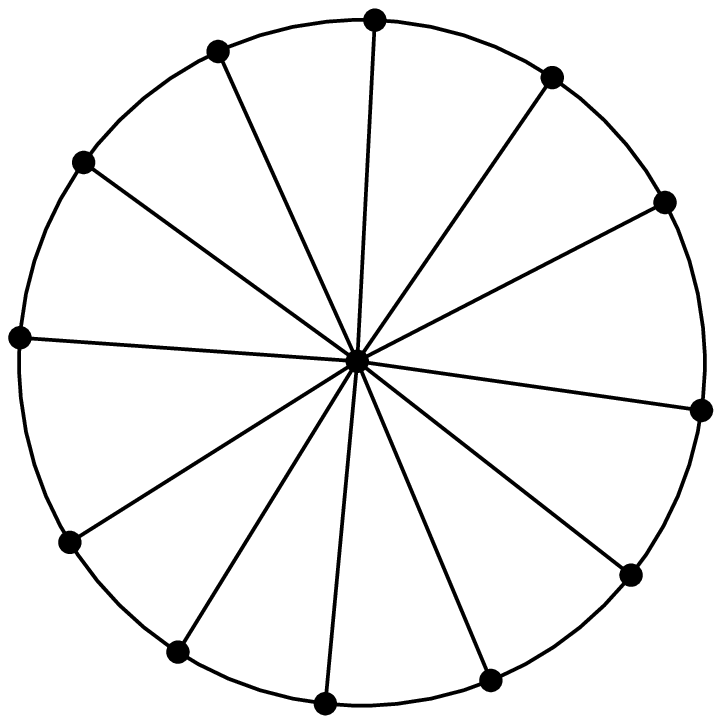}} 
\begin{center}
{\caption{{Wheel graph $ W_{13}=K_1\bigtriangledown C_{12}$ }}
}\end{center}
\end{figure} 
\section{Main Results}

\begin{flushleft}
The aim of this section is to show that any graph cospectral with a multicone graphs $ K_w\bigtriangledown mC_n $  is either regular or bidegreed.
\end{flushleft}

\subsection{Connected graphs cospectral with a multicone graph $ K_w\bigtriangledown mC_n $ with respect to adjacency spectrum.}

\begin{proposition} Let $ G $ be a  graph cospectral with a multicone graph $ K_w\bigtriangledown mC_n$. Then $ {\rm{Spec}}_{A}(G)$ is either
$\bigcup\limits_{k = 1}^{\frac{n}{2} - 1} {\left\{ {{{\left[ { - 1} \right]}^{w-1}},\,\,{{\left[ {{{{2cos}{\dfrac{2k\pi}{n}}}}} \right]}^{2m}}},\,{{\left[ 2 \right]}^{m - 1}},\,{{\left[ { - 2} \right]}^m},\,{{\left[ {\dfrac{{\Omega  - \sqrt {{\Omega ^2} - 4\Gamma } }}{2}} \right]}^1},{{\left[ {\dfrac{{\Omega  + \sqrt {{\Omega ^2} - 4\Gamma } }}{2}} \right]}^1} \right\}}$ or \\
$\bigcup\limits_{k = 1}^{\frac{n-1}{2} } {\left\{ {{{\left[ { - 1} \right]}^{w-1}},\,\,{{\left[ {{{{2cos}{\dfrac{2k\pi}{n}}}}} \right]}^{2m}},\,{{\left[ 2 \right]}^{m - 1}},\,{{\left[ {\dfrac{{\Omega  - \sqrt {{\Omega ^2} - 4\Gamma } }}{2}} \right]}^1},{{\left[ {\dfrac{{\Omega  + \sqrt {{\Omega ^2} - 4\Gamma } }}{2}} \right]}^1}} \right\}}$, where $ \Omega=w+1$ and  $ \Gamma=2(w-1)-mnw $.
\end{proposition}

\noindent \textbf{Proof} It is well-known that the $n$-cycle $C_n$ has eigenvalues $ 2cos\dfrac{2k\pi}{n} $, where $k = 0, . . . , n - 1$. All multiplicities are $ 2 $, except that of $ 2 $ and possibly $ -2 $. Now, by Theorem \ref{the 2-100} the proof is straightforward.
$ \Box $

\begin{lemma}\label{lem 3-1} Let $ G $ be a  graph cospectral with a  multicone graph $ K_w\bigtriangledown mC_n $. Then $ \delta(G) = w+2 $.
\end{lemma}

\noindent \textbf{Proof}
 Let $ \delta(G) = w+2+x$, where $ x $ is an integer number. First, it is clear that in this case the equality in Theorem \ref{the 2-2} happens if and only if $ x=0 $. We claim that $ x=0 $. By contrary, we suppose that $ x\neq 0 $.  Theorem \ref{the 2-2} together with Proposition 3.1 follow  that $ \varrho(G)=\dfrac{w+1+ \sqrt{8k-4l(w+2)+(w+3)^{2} }}{2} \lneqq
\dfrac{w+1+x+ \sqrt{8k-4l(w+2)+(w+3)^{2}+x^{2}+(2w+6-4l)x }}{2} $,  where (as usual) $ k $ and $ l $ denote the numbers of vertices and edges in $ G $, respectively.  

For convenience, we let $ S=8k-4l(w+2)+(w+3)^{2}\geq 0$ and $C = w+3-2l$, and also let $ q(x)=x^{2}+(2w+6-4l)x = x^2 + 2Cx.$ \\ 

Then clearly

\begin{equation*} \sqrt{S}-\sqrt{S+q(x)} < x. \tag{1} \end{equation*}

We consider two cases:

{\bf Case 1.} $x \lneqq 0$.

It is easy and straightforward to see that $ |\sqrt{S}-\sqrt{S+q(x)}| >|x| $, since $ x<0 $.\\
\noindent Transposing and squaring yields

$$2S+q(x)-2\sqrt{S(S+q(x))} > x^{2}.$$

\noindent Replacing $q(x)$ by $x^2 + 2Cx$, we get 

$$S+Cx > \sqrt{S(S+x^2 + 2Cx)}.$$

\noindent Obviously $ Cx\geq 0 $. Squaring again and simplifying yields 

\begin{equation*} C^{2} > S. \tag{2} \end{equation*}

Therefore,

\begin{equation*} k<\frac{l(l-1)}{2}\tag{3}. \end{equation*}

So, if $ x<0 $, then $ G $ cannot be a complete graph. In other words, if $ G $ is a complete graph, then $ x> 0 $. Or one can say that if $ G $ is a complete graph, then $ \delta(G) > w+2 $ ($ \dagger $).

 {\bf Case 2.} $ x\gneqq 0 $.
 
In the same way of {\bf Case 1}, we can conclude that:

If $ G $ is a complete graph, then $ \delta(G) < w+2 $ ($ \dagger\dagger $).

But, two inequalities ($ \dagger $) and ($ \dagger\dagger $) cannot happen together. So we must have $ x=0 $.

 Therefore, the claim holds.
$ \Box $
\begin{lemma}\label{lem 3-2} Let $ G $ be a graph cospectral with a  multicone graph $ K_w\bigtriangledown mC_n $. Then $ G $ is either regular or bidegreed in which any vertex of $ G $ is of degree $ w-1+mn $ or $ w+2 $.
\end{lemma}
\noindent \textbf{Proof}  This follows from Lemma \ref{lem 3-1} together with Theorem \ref{the 2-2}.
$ \Box $

\section{Connected graphs cospectral with the multicone graph $ K_1\bigtriangledown mC_n $. }
\begin{flushleft}
 In this section, we show that any graph cospectral with the  multicone graphs $ K_1\bigtriangledown mC_n $, the cone of graphs $ mC_n $, is isomorphic to $ K_1\bigtriangledown mC_n $. 
\end{flushleft}

\begin{lemma}\label{lem 4-1} Any  graph cospectral with the multicone graph $ K_1\bigtriangledown mC_n $ is DS with respect to its adjacency spectrum.
\end{lemma}
\noindent \textbf{Proof} Let $ G $ be cospectral with the multicone graph $ K_1\bigtriangledown mC_n $. If $ m=1 $ and $ n=3 $ there is nothing to prove, since in this case $ G $ is regular (see Theorem \ref{the 2-5}). Hence we suppose that $ m\neq 1 $ or $ n\neq 3 $. By Lemma \ref{lem 3-2}, it is clear that $ G $ has one vertex of degree $ mn $, say $ j $. We consider
two cases:

{\bf Case 1.} Let $ n $ be even. In this case, it follows from Proposition \ref{prop 2-1}  that 
$ P_{G-j}(x)=(x-\mu_0)^{m-2}\prod\limits_{r=1}^{\frac{n}{2}-1}{{{(x - {\mu _r})}^{2m - 1}}}(x-\mu_{\frac{n}{2}})^{m-1} \sum\limits_{k = 0}^{\frac{{n}}{2}+2} {{\alpha}^2_{kj}{A_k}} $, where $ A_k=\prod\nolimits_{_{i \ne k}i = 0}^{\frac{{n }}{2}+2} {(x - {\mu _i})} $, $ \mu_{\frac{n}{2}}=-2 $, $ \mu_{\frac{n}{2}+1}=1+\sqrt{1+mn} $, $ \mu_{\frac{n}{2}+2}=1-\sqrt{1+mn} $ and  $ \mu _r = 2\cos \frac{{2\pi r }}{n} $ ( $ 0\leq r\leq \frac{n}{2}-1$).\\

 {\bf Case 2.} Let $ n $ be odd. In the same way of {\bf Case 1}, it follows from Proposition \ref{prop 2-1}  that 
 $P_{G-j}(x)={(x - {\mu _0})^{m - 2}}\prod\limits_{r = 1}^{\frac{{n-1}}{2}} {{{(x - {\mu _i})}^{2m - 1}}} \sum\limits_{l = 0}^{\frac{{n -1}}{2}+2} {{\alpha}^2_{lj}{B_l}}$, where $ B_l=\prod\nolimits_{_{i \ne l}i = 0}^{\frac{{n -1}}{2}+2} {(x - {\mu _i})}$, $ \mu_{\frac{n-1}{2}+1}=1+\sqrt{1+mn} $,  $ \mu_{\frac{n-1}{2}+2}=1-\sqrt{1+mn} $ and  $ \mu _r = 2\cos \frac{{2\pi r }}{n} $ ( $ 0\leq r\leq \frac{n-1}{2}$).\\
Now, from Lemma \ref{lem 3-2}, it follows that $ G-j $ is a regular graph and degree of its regularity is $ 2 $. Also, $ G-j $ has $ mn $ vertices of degree $ 2 $. Therefore, we conclude that $ {\rm{Spec}}_{A}(G-j)$ is either $\bigcup\limits_{k = 1}^{\frac{n}{2} - 1} {\left\{ {\,{{\left[ {{{{2cos}{\dfrac{2k\pi}{n}}}}} \right]}^{2m}},\,{{\left[ 2 \right]}^{m },\,{{{\left[ { - 2} \right]}^m}}}} \right\}}$  or $\bigcup\limits_{k = 1}^{\frac{n-1}{2} } {\left\{ {\,{{\left[ {{{{2cos}{\dfrac{2k\pi}{n}}}}} \right]}^{2m}},\,{{\left[ 2 \right]}^{m }}} \right\}}$. Hence $ G-j\cong mC_n $. So, $ G\cong K_1\bigtriangledown mC_n$. 

This completes the proof.
$ \Box $

\begin{flushleft}
Up to now, we have shown that the multicone graph $ K_1\bigtriangledown mC_n $ is DS with respect
to their adjacency spectrum. The natural question is; what happens for  multicone graph $ K_w\bigtriangledown mC_n $? we will respond to this question in the following theorem. 
\end{flushleft}
\begin{theorem}\label{12} Any  graph cospectral with a multicone graph $ K_w\bigtriangledown mC_n $ is isomorphic to $ K_w\bigtriangledown mC_n $.
\end{theorem}
\noindent \textbf{Proof}
We perform the mathematical induction on $  w$. If $ w=1 $, by Lemma \ref{lem 4-1} the proof is clear. We suppose that the claim be valid for $ w $. In  other words, if $ {\rm{Spec}}_{A}(H)= {\rm{Spec}}_{A}(K_w\bigtriangledown mC_n) $, then $ H\cong K_w\bigtriangledown mC_n $, where $ H $ is an arbitrary graph cospectral with a multicone graph $ K_w\bigtriangledown mC_n $. We show that the claim is true for $ w+1 $; that is, we show that if ${\rm{Spec}}_{A}(K)={\rm{Spec}}_{A}(K_{w+1}\bigtriangledown mC_n) $, then $ K\cong K_{w+1}\bigtriangledown mC_n $, where $ K $ denotes a graph cospectral with multicone graph $K_{w+1}\bigtriangledown mC_n$. Graph $ K $ has one vertex and $ w+mn $ edges more than $ H $. By Lemma \ref{lem 3-2}, $ H $ has $ w $ vertices of degree $ w-1+mn $ and $ mn $ vertices of degree $  w+2 $. Also, this lemma implies that $ K $ has $ w+1 $ vertices of degree $ w+mn $ and $ mn $ vertices of degree $  w+3 $. So, we must have $ K\cong K_1\bigtriangledown H $. Now, the inductive hypothesis completes the proof.
$ \Box $
\vspace{3mm}
\begin{flushleft}
In the following, we present an alternative proof of  Theorem \ref{12}.
\end{flushleft}
\noindent \textbf{Proof}(an alternative proof of Theorem \ref{12}) Let $ G $ be a graph cospectral with a multicone graph $ K_w\bigtriangledown mC_n $.  By Lemma \ref{lem 3-2}, $ G $ consists of subgraph $ \Gamma $ and degree of any  vertex of $ \Gamma $ is $ w-1+mn $. In  other words, $ G\cong K_w\bigtriangledown H $, where $ H $ is a subgraph of $ G $. Now, we remove  vertices of the complete graph $ K_w $ and we  consider $ mn $ another vertices. Consider  subgraph $ H $ consisting of these $ mn $ vertices. $H$ is regular and  degree of its regularity is $ 2 $ and multiplicity of $ 2 $ is $ m $. In other words, $ H $ is a cycle or it is a disjoint union of  several $n$-cycles. By Theorem \ref{the 2-100},  $ {\rm{Spec}}_{A}(H)$ is either  $\bigcup\limits_{k = 1}^{\frac{n}{2} - 1} {\left\{ {\,{{\left[ {{{{2cos}{\dfrac{2k\pi}{n}}}}} \right]}^{2m}},\,{{\left[ 2 \right]}^{m}}\,{{\left[ { - 2} \right]}^m}} \right\}}$  or $\bigcup\limits_{k = 1}^{\frac{n-1}{2} } {\left\{ {\,{{\left[ {{{{2cos}{\dfrac{2k\pi}{n}}}}}\right]}^{2m}}},\,{{\left[ 2 \right]}^{m }} \right\}}$. Hence $ {\rm{Spec}}_{A}(H)={\rm{Spec}}_{A}(mC_n) $.  
This follows the result.
$ \Box $

\section{  Some classes of graphs $\overline{K_w\bigtriangledown mC_n} $ and  their adjacency spectrum}
\begin{flushleft}
In this section, we show that graphs $ \overline{K_w\bigtriangledown mC_3} $ are DS with respect to their adjacency spectra.
\end{flushleft}
\begin{theorem} Let $ G $ be a graph and $ {\rm{Spec}}_{A}(G)={\rm{Spec}}_{A}(\overline{K_w\bigtriangledown mC_3}) $. Then $ G\cong \overline{K_w\bigtriangledown mC_3} $.
\end{theorem}

\noindent \textbf{Proof}
It is clear that $ {\rm{Spec}}_{A}(G)=\left\{ {{{\left[ { - 3} \right]}^{m - 1}},\,{{\left[ 0 \right]}^{2m + w}},\,{{\left[ {3m - 3} \right]}^1}} \right\}$. If $ m=1 $, there is nothing to prove. If $ m=2 $, Theorem \ref{the 2-7} implies that $ G\cong G_1\cup sK_1 $, where $ s $ and $ G_1 $ are  a natural number and a complete multipartite graph, respectively. On the other hand, it follows from Theorem \ref{the 2-6} that $ G_1 $  is  a complete bipartite graph.  Therefore, $ G_1\cong K_{3, 3}$. Hence $ G\cong wK_1 \cup K_{3, 3} $. Therefore, we can suppose that $ m\geqslant 3 $. We know that regularity and the number of triangles of graph $ G $ can be determined by its adjacency spectrum. So, $ G $ is a non-regular graph and it is not a complete bipartite graph. Therefore, from Theorems \ref{the 2-6} and \ref{the 2-7} we conclude that $ G $ is not connected and $ G=G_1\cup mK_1 $, where $ m $ is a natural number. Also, it is clear that $ G_1 $ has exactly three distinct eigenvalues. Now,  by contrary, we suppose that $ G_1 $ be a non-regular graph. In this case, it follows from Theorem \ref{the 2-6}  that $ G_1 $ is a complete bipartite graph and so $ G $ is a bipartite graph. This is a contradiction. Hence $ G_1 $ is a regular graph and so $ \varrho(G_1)=3m-3 $. It  follows from Theorem \ref{the 2-7}  that $ G_1 $ is a multipartite graph. Hence $ G_1\cong K_{\underbrace { 3, 3, ..., 3}_{m\,times}}$ and so $ G=G_1 \cup wK_1\cong K_{\underbrace { 3, 3, ..., 3}_{m\,times}}\cup wK_1\cong\overline {K_w \bigtriangledown mK_3}$. 
This follows the result.
$ \Box $

\section{ Connected graphs cospectral with a multicone graph $ K_w\bigtriangledown mC_n$ with respect to Laplacian spectrum}
\begin{flushleft}
In this section, we show that if $ w, m\neq 1 $ and $ n\neq 6 $ (see 19, Fig. 1), then multicone graphs $ K_w\bigtriangledown mC_n$ are DS with respect to their Laplacian spectrum.
\end{flushleft}
\begin{theorem}\label{14}
Multicone graphs $ K_w\bigtriangledown mC_n $ are DS with respect to their Laplacian spectrum, where $ w, m\neq 1 $ and $ n\neq 6 $.
\end{theorem}

\noindent \textbf{Proof} It is well-known that ${\rm{Spec}}_{L}(C_n)= 2-2cos\dfrac{2k\pi}{n} $, where $k = 0, . . . , n - 1$. All multiplicities are $ 2 $, except that of $ 0 $ and possibly $ 4$. We perform mathematical induction on $  w$. First, we suppose that $ n $ be even. If $ w=1 $, by Proposition \ref{prop 2-2} the proof is straightforward. Let the claim be true for $ w $, that is, $ {\rm{Spec}}_{L}(G_1)= {\rm{Spec}}_{L}(K_w\bigtriangledown mC_n)= 
\bigcup\limits_{k = 1}^{\frac{n}{2} - 1} {\left\{ {\,{{\left[ {w + mn} \right]}^w}\,,{{\left[ {w + 2 - 2\cos \dfrac{{2\pi k}}{n}} \right]}^{2m}},\,{{\left[ w \right]}^{m - 1}},\,{{\left[ {w + 4} \right]}^m},\,{{\left[ 0 \right]}^1}} \right\}}$ follows that $ G_1\cong K_w\bigtriangledown mC_n $. We show that the claim is true for $ w+1 $, that is, we show that $ {\rm{Spec}}_{L}(G)={\rm{Spec}}_{L}(K_{w+1}\bigtriangledown mC_n)=
\bigcup \limits_{k = 1}^{\frac{n}{2} - 1} {\left\{ {\,{{\left[ {w +1+ mn} \right]}^{w+1}}\,,{{\left[ {w + 3 - 2\cos \dfrac{{2\pi k}}{n}} \right]}^{2m}},\,{{\left[ w+1 \right]}^{m - 1}},\,{{\left[ {w + 5} \right]}^m},\,{{\left[ 0 \right]}^1}} \right\}} $ follows that $ G\cong K_{w+1}\bigtriangledown mC_n$. It is clear that, $ G $ has one vertex and $w+ mn$
edges more than $ G_1 $. On the other hand, Theorem \ref{the 2-4}  follows that any of graphs $ G $ and $ G_1 $ is the join of two graphs. In addition, ${\rm{Spec}}_{L}(K_1\bigtriangledown G_1)={\rm{Spec}}_{L}(G)  $. Therefore, we must have $ G\cong K_1\bigtriangledown G_1 $. Now, the induction hypothsis completes the proof. If $ n $ is odd, the proof is similar.
$ \Box $

\section{Some algebraic properties about multicone graphs  $ K_w\bigtriangledown mC_n $.}
\begin{flushleft}
It is proved that a graph $ G $ is perfect if and only if $ G $ is Berge; that is, it contains no odd hole or antihole as induced subgraph, where odd hole and antihole are odd cycle, $ C_m $ for $ m\ge 5 $, and its complement, respectively. Also, in $1972$ Lov\'asz proved that, a graph is perfect if and only if its complement is perfect (see \cite{B123}).\\ In this  section, we show that if $ n $ is either even or  $ 3 $, then any graph cospectral with $K_w\bigtriangledown mC_n $ with respect to its adjacency spectrum as well as its Laplacian spectrum must be perfect. 
\end{flushleft}
\begin{theorem}\label{140} Let graph $ G $ be cospectral with a multicone graph $ K_w\bigtriangledown mC_n $. Then:

\begin{flushleft}
 $ G $ and $ \overline{G} $ are perfect if and only if $ n $ is either even or  $ 3 $. 
\end{flushleft}
\end{theorem}
\noindent \textbf{Proof}
$(\Rightarrow)$ By what were said  in the beginning of this section and Theorem \ref{12} the proof is straightforward.

$( \Leftarrow) $ It is quite clear that $ G $ cannot consist of an odd hole of order greater than or equal to five as an induced subgraph. We show that $ G $ contains no odd antihole of order greater than or equal to five as an induced subgraph. By contrary, we suppose that $ G $ contains $ \overline{C_k}$ as an induced subgraph, where $ k $ is an odd natural number greater than or equal to five. Hence $ \overline{G}=wK_1\cup \overline{mC_n} $ must consists of $ C_k $ as an induced subgraph. In other words, $  \overline{mC_n}=\underbrace {\overline{C_n}\bigtriangledown .\,.\,.\bigtriangledown \overline{C_n}}_{m\,\,times}$ must consists of $ C_k $ as an induced subgraph.  This is obviously a contradiction.

This observation completes the proof.
$ \Box $

\begin{theorem}Let $ G $ be a graph and  $ {\rm{Spec}}_{L}(G)={\rm{Spec}}_{L}(K_w\bigtriangledown mC_n) $. Then:

\begin{flushleft}
 $ G $ and $ \overline{G} $ are perfect if and only if $ n $ is either even or  $ 3 $.
\end{flushleft}
\end{theorem}
\noindent \textbf{Proof}
The proof is in a similar manner of Theorem \ref{140}.
$ \Box $

\section{ Two Conjectures}
In this paper, it were proved that connected multicone graphs $ K_w\bigtriangledown mC_n $ are determined by both their adjacency and Laplacian spectra. Also, we show that $ \overline{ K_w\bigtriangledown mC_3} $ are determined by their adjacency spectra. By \cite{A} (Theorem 5.2), we can deduce that graphs $ \overline{ K_w\bigtriangledown C_4} $ are determined by their adjacency spectra. Now, we pose the following conjectures.

\begin{conjecture}  Graphs $\overline{ K_w\bigtriangledown mC_n} $ are DS with respect to their adjacency spectrum.
\end{conjecture}
\begin{conjecture}  Multicone graphs $ K_w\bigtriangledown mC_n $ are DS with respect to their signless Laplacian spectrum.

\end{conjecture}

\end{document}